\documentclass[11pt]{article}

\usepackage[english]{babel}
\usepackage{graphicx}
\usepackage{subcaption}

\newcommand{\R}{\mathbb{R}}

\newcommand{\C}{\mathbb{C}}
\newcommand*{\bb}[1]{\mathbf{#1}}

\usepackage{siunitx}

\usepackage[margin=2.5cm]{geometry}
 \usepackage{amsmath, amssymb}
\usepackage{amsthm}

\usepackage[table]{xcolor}
\usepackage[colorlinks = true,
            linkcolor = blue,
            urlcolor  = blue,
            citecolor = blue,
            anchorcolor = blue]{hyperref}
      
\newtheorem{theorem}{Theorem}[section]
\newtheorem{remark}[theorem]{Remark}

\usepackage{titling}

\usepackage{array,multirow}
\usepackage{hhline}
{

\begin{document}

\title{Rainfall infiltration: Direct and Inverse problems on a linear evolution equation}

\author{Konstantinos Kalimeris\thanks{Mathematics Research Center, Academy of Athens, Greece} \, and Leonidas  Mindrinos\thanks
  {Department of Natural Resources Development and Agricultural Engineering, Agricultural University of Athens, Greece}}

\thanksmarkseries{arabic}

\maketitle

\begin{abstract}
Originating from the mathematical modelling of rainfall infiltration, we derive the solution of an initial-boundary value problem of a linear evolution partial differential equation, by using the Fokas method. We present numerical examples which correspond to specific physical rainfall problems. Based on this formalism we present an effective algorithm for the associated null-controllability problem, namely
 we numerically derive a family of boundary controls that steer the solution to the desired flat final state. Finally, a regularisation scheme allows the derivation of relatively small controls, in cases where this is necessary.
\end{abstract}

\section{Introduction and Mathematical Formulation}\label{intro}

A plethora of physical phenomena is mathematically modelled as initial and boundary value problems (IBVPs) on linear and nonlinear evolution partial differential equations (PDEs). The interest of providing analytical solutions on such problems extends also to the fact that this provides a background to numerically solve inverse problems associated to the former ones. Usually,  linear problems yield  simple (or even simplistic) modelling of the detailed description of a physical phenomenon. However, most of the times this (linear) step is indispensable in the comprehension of the mathematical approach which will effectively address the original problem using both analytical and numerical approaches.

Modelling soil infiltration, which describes  the process of water flow into the soil through its surface, induces  a large class of such problems; for more details on this process, as well as how this is modelled to an IBVP for an evolution PDE, we direct the reader to \cite{arg24,Bra73,hayek16, Philip91,Tracy11}. Therein, one could see the details of obtaining  \eqref{bvp1} as the main linear PDE which describes a simple model of vertical infiltration in bounded homogeneous profiles. In \cite{arg24} the association of different classes of initial-boundary conditions to specific physical phenomena is provided. The motivation of the current work, namely considering the IBVP \eqref{bvp} is given in the paragraph below equation \eqref{bvp}.

One other important aspect associated to evolution problems is the attempt to control the relevant procedure. This is a typical situation of an inverse problem which has attracted a lot of interest both analytically and numerically. It is not in the scope of this work to provide a concise overview of the different set-ups of controllability, rather we focus our effort to the so-called boundary null controllability. This problem consists of constructing appropriate boundary conditions which will steer the solution of the IBVP to a zero-state, for any given initial condition: of course, this is not always doable, see for example the lack of null controllability of the heat equation on the half-line \cite{zua01}. To the best of our knowledge, in the case that a numerical solution can be provided there are three avenues to address the problem: (i) an indirect one, based on the Hilbert Uniqueness method \cite{Ros98}, (ii) a direct one, based on the `flatness aproach' and the smoothness of the initial data \cite{MarRos14} (iii) a direct one, based on the Fokas method \cite{KOD23}. Considering an overview of the implementation of these methods, as well as a brief discussion on their advantages and disadvantages we refer to \cite{KOD23}. In fact, this work introduced for the first time the implementation of the Fokas method to the numerical construction of  boundary controls. It is a part of a trilogy by \"Ozsar\i   \ and one of the authors which introduces the implementation of the Fokas method to the area of boundary control problems. In \cite{KO20} an elementary proof was provided for the well known result of lack of null controllability of the heat equation on the half-line \cite{zua01}, which is also much shorter than the one existed in the literature at that time. In \cite{OK23} an analytical proof on the open question of controllability of two  Schr\"odinger-type equations was provided.

In the current work, we are extending the methodology of \cite{KOD23} to null controlling \eqref{bvp}. We mention that this kind of extension has already been performed in \cite{hwang23}, but for simpler boundary conditions and controls. Nevertheless, we point out the importance of the current work for three reasons: First, the physical problem that we present and solve here cannot be captured by the setup in  \cite{hwang23}. Second, we would like to clarify some potential misleads on the performance of the algorithm (namely, the exponential decay of the error in terms of the number of discretization  points), as well as the possibility of controlling the IBVP in very short time. Third, we introduce a regularisation technique which allows the construction of relatively small controls. In contrast with what is vaguely implied in \cite{hwang23}, this IBVP can be controlled in arbitrarily small time, $T>0$, for any value of the water diffusivity parameter $D_0>0$, with the cost that the control is growing (sub)exponentially as $T$ tends to 0. For an analytical derivation of this result we refer to survey \cite{zua07}, the recent works \cite{galo24a, galo24b} and the literature therein, where $O\left(e^{-c/T}\right)$ is proven as a bound for the cost of the control, for some constant $c>0$ which depends on the other parameters of the IBVP; we note that preliminary numerical results in \cite{KOD23} indicate a slightly better conjecture, namely $O\left(e^{-c/\sqrt{T}}\right)$, which is compatible with the results and the discussion in \cite{galo24b}. Furthermore, in the current work we are able to mitigate this effect by loosening the demand of strictly zero deviation of the final solution.

The current work is based on the Fokas method for solving linear IBVPs. This method, which is also known as the unified transform, was introduced in 1997 by Fokas \cite{F97}, for solving IBVPs for nonlinear integrable PDEs. Later, it was realised that it produced effective analytical and numerical solutions for linear PDEs \cite{F02}; for an overview of the method we refer to \cite{F08, FK22}. A very interesting outcome/discussion on the importance of the solutions of the linear problems to non-linear ones can be found on works associated to regularity results \cite{BFO20, Him20, OY19}. The analytical solution of linear problems provided by the Fokas method is in integral form, also for bounded domains (in the current work this is a finite interval), which involve exponentially decaying integrands along the contours of integration. Using this feature it is obvious that the numerical computation of the solution becomes very efficient, as opposed to existing solutions which are in the form of oscillatory series. For the same reason, the control problem that we study here,  displays a straightforward formulation and rather effective implementation. One could identify the analytical source of these numerical feature in the uniform convergence of the solution (of the direct problem) at the boundary of the domain: one can retrieve the boundary values of the relevant IBVP via a direct calculation of the integral representation of the solution at the boundary (in the current work at the endpoints of the interval), which is to be contrasted with classical  representations of the solution.

We are interested in solving the following initial boundary value problem
\begin{subequations}\label{bvp}
\begin{alignat}{3}
\frac{\partial \theta}{\partial t} + K_0 \frac{\partial \theta}{\partial x} &= D_0 \frac{\partial^2 \theta}{\partial x^2},  \quad && 0<x <L, \, t>0,  \label{bvp1}\\
\theta (x,0) &= \theta_0 (x), \quad && 0 <x <L, \label{bvp2}\\ 
\theta (0,t) -\alpha \frac{\partial \theta}{\partial x} (0,t) &= f(t),  \quad &&t>0, \label{bvp3}\\
\theta (L,t) -\beta \frac{\partial \theta}{\partial x} (L,t) &= g(t),  \quad &&t>0, \label{bvp4}
\end{alignat}
\end{subequations}
for the water content $\theta$ in a bounded and homogeneous soil with length $L>0$ characterized by the water diffusivity $D_0$ and the hydraulic conductivity $K_0.$ The initial soil water content $\theta_0$, the coefficients $\alpha$ and $
\beta$ and the time-dependent boundary functions $f$ and $g$ are arbitrary but given.

The motivation to consider the above IBVP with Robin-Robin boundary conditions comes from its appearance in the mathematical modeling of vertical infiltration. The boundary conditions for $\alpha = \beta = \tfrac{D_0}{K_0}$ model rainfall infiltration. If the boundary function ($f$ or/and $g$) is time-independent we consider rainfall with constant flux and if it is zero we describe an impermeable surface either at $x=0$ or at $x=L$ \cite{Bra73, Philip91}. Here, we do not fix $\alpha$ and $\beta$ from the beginning in order to examine a more general problem that can be reduced to a Robin-Dirichlet problem or a Dirichlet-Dirichlet problem as in \cite{arg24}. Thus, the solution of the direct problem gives the evolution of the distribution of the water content in the soil during rainfall.

The inverse problem which we consider in this work involves the boundary control of this IBVP: find the appropriate boundary condition which will steer the solution to the desired flat state. The physical meaning involves the regulation of the permeability of the bottom surface which would lead to the desired distribution of the water content in the soil, in the event of a rainfall. A simple change of variables in the IBVP (resp. in the solution) can be interpreted as an adjustment to the rainfall and/or soil conditions on the top surface, leading to the desired distribution.

\subsection*{Organisation of the paper}
In the next section we solve analytically the direct problem, namely we provide the solution of the IBPV \eqref{bvp} in the form of an integral representation of (transforms of) the initial and boundary data; we also provide numerical examples which come from specific physical applications. In section \ref{sec_inv}, we present and implement the numerical algorithm for constructing the control,  for different values of the parameters of the problem. In more detail, the reader can  find numerical evidence for the exponential decay of the numerical error as the discretization points increase. Finally, one can observe the increase that the control suffers from, as the diffusivity, $D_0>0$, or/and the duration of the action of the control, $T>0$, become small, as well as the mitigation of this effect by the regularised controls.

\section{Direct problem}\label{sec_formulation}

\subsection{Analytical solution}

We apply the Fokas method to solve \eqref{bvp}, which  consists of three steps \cite{F08,FK22} (i) construct the `global relation' of the initial and boundary data with unknown boundary conditions, (ii) derive an `integral representation' for the solution which involves both known and unknown initial and boundary values, (iii) employ the symmetries of the global relation to evaluate the contribution of the unknown boundary values in the integral representation. This procedure will provide the solution of \eqref{bvp} as integral representation of (tranforms of) the  initial and boundary data.

Indeed, we multiple \eqref{bvp1} with $e^{-i \lambda x}$ and we integrate with respect to $x,$ to obtain
\begin{equation}\label{integral}
\frac{\partial}{\partial t} \int_0^L e^{-i \lambda x} \theta (x,t) dx + K_0 \int_0^L e^{-i \lambda x} \frac{\partial \theta}{\partial x} (x,t) dx  = D_0 \int_0^L e^{-i \lambda x} \frac{\partial^2 \theta}{\partial x^2} (x,t) dx.
\end{equation}
We define 
\begin{equation}\label{fourier}
\hat\theta (\lambda,t) =  \int_0^L e^{-i \lambda x} \theta (x,t) dx
\end{equation}
and perform integration by parts to the other two terms to obtain
\begin{equation}\label{parts1}
\int_0^L e^{-i \lambda x} \frac{\partial \theta}{\partial x} (x,t) dx = e^{-i \lambda L} \theta (L,t) - \theta (0,t) + i \lambda \hat\theta (\lambda,t)
\end{equation}
and
\begin{equation}\label{parts2}
\int_0^L e^{-i \lambda x} \frac{\partial^2 \theta}{\partial x^2} (x,t) dx = e^{-i \lambda L}
\left(\frac{\partial \theta}{\partial x} (L,t) +i\lambda  \theta (L,t) \right) - \frac{\partial \theta}{\partial x} (0,t) - i\lambda  \theta (0,t) -\lambda^2 \hat\theta (\lambda,t).
\end{equation}
We substitute \eqref{parts1} and \eqref{parts2} in \eqref{integral} to get
\begin{equation}\label{integral2}
\frac{\partial}{\partial t}\hat\theta (\lambda,t) + \omega (\lambda) \hat\theta (\lambda,t) = D_0 e^{-i \lambda L}  \left( \frac{\partial \theta}{\partial x} (L,t)- i \nu \theta (L,t)\right) - D_0  \left( \frac{\partial \theta}{\partial x} (0,t)- i \nu \theta (0,t)\right),
\end{equation}
where 
\begin{equation*}
\omega (\lambda) =  D_0  \lambda^2 + i K_0 \lambda \quad \text{ and } \quad \nu =\nu(\lambda) = -\lambda-i\tfrac{K_0}{D_0}.
\end{equation*}

We define
\begin{equation}\label{t_transform}
\begin{aligned}
\tilde f_0 (\lambda,t) &= \int_0^t e^{\lambda \tau} \theta (0,\tau ) d\tau,\quad
\tilde f_1 (\lambda,t) = \int_0^t e^{\lambda \tau} \frac{\partial \theta}{\partial x} (0,\tau )d\tau, \\
\tilde g_0 (\lambda,t) &= \int_0^t e^{\lambda \tau} \theta (L,\tau ) d\tau, \quad
\tilde g_1 (\lambda,t) = \int_0^t e^{\lambda \tau} \frac{\partial \theta}{\partial x} (L,\tau )d\tau,
\end{aligned}
\end{equation}
the so-called $t-$transform of the boundary values and we multiple \eqref{integral2} with $e^{\omega(\lambda)t}$ and we integrate with respect to $t,$ to derive
\begin{equation}\label{global0}
\begin{aligned}
\hat\theta (\lambda,t)  &= e^{-\omega (\lambda) t} \left[
  \hat \theta_0 (\lambda) + D_0 e^{-i \lambda L} \left( \tilde g_1 (\omega (\lambda),t) - i \nu \tilde g_0 (\omega (\lambda),t)  \right) \right.\\
&\phantom{=}\left.   -  D_0 \left( \tilde f_1 (\omega (\lambda),t) - i \nu \tilde f_0 (\omega (\lambda),t)  \right) \right],
  \end{aligned}
\end{equation}
where $\hat \theta_0 (\lambda) := \hat \theta (\lambda,0).$ In order to  consider the Robin boundary conditions we first apply the $t-$transform in \eqref{bvp3} and \eqref{bvp4} to obtain
\begin{equation}\label{t_robin}
\tilde f_1  = \frac{1}{\alpha} ( \tilde f_0  -  \tilde f), \quad \tilde g_1  = \frac{1}{\beta} (  \tilde g_0  -  \tilde g),
\end{equation}
where $\tilde{f}$ and $\tilde{g}$ denote the $t-$transforms of the boundary functions $f$ and $g,$ respectively.
Then, the substitution of \eqref{t_robin} in \eqref{global0} results in
\begin{equation}\label{global}
\begin{aligned}
\hat\theta (\lambda,t)  &= e^{-\omega (\lambda) t} \left\{
  \hat \theta_0 (\lambda) + D_0 e^{-i \lambda L} \left[ \left( \frac{1}{\beta} - i \nu  \right) \tilde g_0 (\omega (\lambda),t) - \frac{1}{\beta} \tilde g (\omega (\lambda),t) \right] \right.\\
&\phantom{=}\left.   -  D_0 \left[ \left( \frac{1}{\alpha} - i \nu  \right) \tilde f_0 (\omega (\lambda),t) - \frac{1}{\alpha} \tilde f (\omega (\lambda),t) \right] \right\}.
  \end{aligned}
\end{equation}
 This is the global relation connecting the initial and boundary values. The integral representation of the solution follows from the application of the inverse Fourier transform in \eqref{global} \cite{FK22,arg24,barros19}
\begin{equation}\label{repre1}
\begin{aligned}
\theta (x,t) &= \frac1{2\pi} \int_{\R} e^{i \lambda x - \omega (\lambda)t} \hat \theta_0 (\lambda) d\lambda \\
&\phantom{=}- \frac1{2\pi} \int_{\partial D_+} e^{i \lambda x - \omega (\lambda)t} \left[ D_0\left( \frac{1}{\alpha} - i \nu  \right) \tilde f_0 (\omega (\lambda),t) - \frac{D_0}{\alpha} \tilde f (\omega (\lambda),t) \right]d\lambda \\
&\phantom{=}- \frac1{2\pi} \int_{\partial D_-} e^{-i \lambda (L-x) - \omega (\lambda)t} \left[ D_0\left( \frac{1}{\beta} - i \nu  \right) \tilde g_0 (\omega (\lambda),t) - \frac{D_0}{\beta} \tilde g (\omega (\lambda),t)
 \right]d\lambda,
\end{aligned}
\end{equation}
where $\partial D_{\pm} = \big\{ \lambda \in \C: \operatorname{Re}(\omega (\lambda)) = 0, \ \operatorname{Im}(\lambda) \gtrless 0 \big\}$ is the boundary of the domain $D_{\pm} = \big\{ \lambda \in \C: \operatorname{Re}(\omega (\lambda)) < 0, \ \operatorname{Im}(\lambda) \gtrless 0 \big\}.$

\begin{remark}\label{ref-deform}
In fact, there is plenty of freedom on the choice of $\partial D_{\pm}$, being any simply connected open curve which satisfies the following two properties:
\begin{itemize}
\item[a)] they belong to $\C_\pm:= \big\{ \lambda \in \C: \operatorname{Im}(\lambda) \gtrless 0 \big\}$, respectively.
\item[b)] their (linear) asymptotes belong in $\C_\pm / D_{\pm}$, respectively. In particular, the asymptotes are defined by $A_{\pm} = \big\{ \lambda \in \C: \operatorname{Im}(\lambda)= c \operatorname{Re}(\lambda), \ c\in[-1,1],\ \ \operatorname{Im}(\lambda) \gtrless 0  \big\}$.
\end{itemize}
This freedom of choice is useful if one would wish to avoid the computation of any pole contribution occurring in the process described in the remaining part of this section. In what follows, we allow ourselves the abuse of notation where with $\partial D_{\pm}$ we denote any favourable such curve; not strictly the boundary of $D_{\pm}$.
\end{remark}

In the representation \eqref{repre1} the two terms $\tilde f_0$ and $\tilde g_0$ are unknown but we can eliminate them using the invariant transformation $\lambda \mapsto \nu =  \nu(\lambda)$ in the global relation \eqref{global}. Given that $\omega (\lambda) =\omega (\nu),$ we obtain
\begin{equation}\label{global2}
\begin{aligned}
\hat\theta (\nu,t)  &= e^{-\omega (\lambda) t} \left\{
  \hat \theta_0 (\nu) + D_0 e^{-i \nu L} \left[ \left( \frac{1}{\beta} - i \lambda  \right) \tilde g_0 (\omega (\lambda),t) - \frac{1}{\beta} \tilde g (\omega (\lambda),t) \right] \right.\\
&\phantom{=}\left.   -  D_0 \left[ \left( \frac{1}{\alpha} - i \lambda  \right) \tilde f_0 (\omega (\lambda),t) - \frac{1}{\alpha} \tilde f (\omega (\lambda),t) \right] \right\}.
  \end{aligned}
\end{equation}

The functions $\tilde f_0$ and $\tilde g_0$ are the solutions of the system  \eqref{global} and \eqref{global2}. We define
\begin{equation}\label{def:Delta}
\Delta (\lambda) = e^{-i \lambda L} (1- i \alpha\lambda) (1-i\beta \nu) - e^{-i \nu L} (1- i \alpha\nu) (1-i\beta \lambda)
\end{equation}
and we derive the formulas: 
\begin{equation}\label{sol_system}
\begin{aligned}
\tilde g_0 &= \frac{1}{\Delta (\lambda)} \left\{  \frac{\beta}{D_0} e^{\omega(\lambda)t} \left[ (1- i \alpha\lambda) \hat \theta (\lambda,t) - (1- i \alpha\nu)\hat \theta (\nu,t)   \right] -\frac{\beta}{D_0} \left[ (1- i \alpha\lambda) \hat \theta_0 (\lambda) - (1- i \alpha\nu)\hat \theta_0 (\nu)   \right] \right. \\
&\phantom{=}\left. +\left[ e^{-i \lambda L} (1- i \alpha\lambda) - e^{-i \nu L} (1- i \alpha\nu)\right] \tilde g (\omega (\lambda),t) -\beta \left(\frac{K_0}{D_0} - i 2\lambda \right) \tilde f (\omega (\lambda),t)
\right\}, \\
\tilde f_0 &= \frac{1}{\Delta (\lambda)} \left\{  \frac{\alpha}{D_0} e^{\omega(\lambda)t} \left[ e^{-i\nu L}(1- i \beta\lambda) \hat \theta (\lambda,t) -e^{-i\lambda L} (1- i \beta\nu)\hat \theta (\nu,t)   \right] \right. \\
&\phantom{=}\left.-\frac{\alpha}{D_0} \left[ e^{-i\nu L}(1- i \beta\lambda) \hat \theta_0 (\lambda) - e^{-i\lambda L}(1- i \beta\nu)\hat \theta_0 (\nu)   \right] \right. \\
&\phantom{=}\left. +\alpha e^{-\tfrac{K_0}{D_0}L}\left( \frac{K_0}{D_0} - i 2\lambda \right) \tilde g (\omega (\lambda),t) - \left[ e^{-i \nu L} (1- i \beta\lambda) - e^{-i \lambda L} (1- i \beta\nu)\right] \tilde f (\omega (\lambda),t)
\right\}.
\end{aligned}
\end{equation}
Substituting \eqref{sol_system} in \eqref{repre1} we get 
\begin{equation}\label{repre2}
\begin{aligned}
\theta (x,t) &= \frac1{2\pi} \int_{\R} e^{i \lambda x - \omega (\lambda)t} \hat \theta_0 (\lambda) d\lambda \\
&\phantom{=}- \frac1{2\pi} \int_{\partial D_+} \frac{e^{i \lambda x}}{\Delta (\lambda)}  \left[ e^{-i\nu L}(1- i \beta\lambda)(1- i \alpha\nu) \hat \theta (\lambda,t) -e^{-i\lambda L} (1- i \beta\nu)(1- i \alpha\nu)\hat \theta (\nu,t)   \right] d\lambda
\\
&\phantom{=}+ \frac1{2\pi} \int_{\partial D_+} \frac{e^{i \lambda x - \omega(\lambda)t}}{\Delta (\lambda)}  \left[ e^{-i\nu L}(1- i \beta\lambda)(1- i \alpha\nu) \hat \theta_0 (\lambda) -e^{-i\lambda L} (1- i \beta\nu)(1- i \alpha\nu)\hat \theta_0 (\nu)  \right] d\lambda \\
&\phantom{=}- \frac{D_0}{2\pi} e^{-\tfrac{K_0}{D_0}L}\int_{\partial D_+} \frac{e^{i \lambda x - \omega(\lambda)t}}{\Delta (\lambda)}(1- i \alpha\nu) \left( \frac{K_0}{D_0} - i 2\lambda \right) \tilde g (\omega (\lambda),t) d\lambda 
\\
&\phantom{=}+ \frac{D_0}{2\pi} \int_{\partial D_+} \frac{e^{-i \lambda (L-x) - \omega(\lambda)t}}{\Delta (\lambda)}(1- i \beta\nu)\left( \frac{K_0}{D_0} - i 2\lambda \right) \tilde f (\omega (\lambda),t) d\lambda 
\\
&\phantom{=} - \frac1{2\pi} \int_{\partial D_-} \frac{e^{-i \lambda (L-x) }}{\Delta (\lambda)}
\left[  (1- i \beta\nu)(1- i \alpha\lambda) \hat \theta (\lambda,t) - (1- i \beta\nu)(1- i \alpha\nu)\hat \theta (\nu,t) 
\right] d\lambda \\
&\phantom{=} + \frac1{2\pi} \int_{\partial D_-} \frac{e^{-i \lambda (L-x) -\omega (\lambda)t }}{\Delta (\lambda)}
\left[  (1- i \beta\nu)(1- i \alpha\lambda) \hat \theta_0 (\lambda) - (1- i \beta\nu)(1- i \alpha\nu)\hat \theta_0 (\nu) 
\right] d\lambda \\
&\phantom{=} - \frac{D_0}{2\pi} e^{-\tfrac{K_0}{D_0}L} \int_{\partial D_-} \frac{e^{i \lambda x -\omega (\lambda)t }}{\Delta (\lambda)}  (1- i \alpha\nu)\left( \frac{K_0}{D_0} - i 2\lambda \right) \tilde g (\omega (\lambda),t) d\lambda  \\
&\phantom{=} + \frac{D_0}{2\pi} \int_{\partial D_-} \frac{e^{-i \lambda (L-x) -\omega (\lambda)t }}{\Delta (\lambda)} (1- i \beta\nu) \left( \frac{K_0}{D_0} - i 2\lambda \right) \tilde f  (\omega (\lambda),t) d\lambda.
\end{aligned} 
\end{equation}
All the integrals of the RHS of this expression involve transforms of known initial and boundary values, apart from the second and sixth integrals, namely the ones that involve the terms $\hat \theta (\lambda,t) $ and $\hat \theta (\nu,t) $. We now show that the latter integrals vanish. Indeed, we examine the contribution of the second integral in the RHS of the above representation. Recall that on $\partial D_+$ since $\operatorname{Im}(\lambda) >0,$ the term $e^{i\lambda L}$ is bounded and vanishes exponentially as $|\lambda| \rightarrow \infty.$ Thus, we obtain the following asymptotic behaviour:
\begin{equation*}
\begin{aligned}
\frac{1}{\Delta (\lambda)}  &\left[ e^{-i\nu L}(1- i \beta\lambda)(1- i \alpha\nu) \hat \theta (\lambda,t) -e^{-i\lambda L} (1- i \beta\nu)(1- i \alpha\nu)\hat \theta (\nu,t)   \right] \\
&\sim  e^{-i(\nu - \lambda)L} \frac{(1- i \beta\lambda)(1- i \alpha\nu)}{(1- i \beta\nu)(1- i \alpha\lambda)} \hat \theta (\lambda,t) -\frac{(1- i \alpha\nu)}{(1- i \alpha\lambda)}\hat \theta (\nu,t)  \\
&= e^{i \lambda L -\tfrac{K_0}{D_0}L} \frac{(1- i \beta\lambda)(1- i \alpha\nu)}{(1- i \beta\nu)(1- i \alpha\lambda)} \int_0^L e^{-i\lambda (\xi - L)} \theta(\xi,t ) d\xi -\frac{(1- i \alpha\nu)}{(1- i \alpha\lambda)}\hat \theta (\nu,t).
\end{aligned} 
\end{equation*}
All terms and the integral, since $\xi-L<0,$ are bounded, namely $O\left(\frac{1}{\lambda}\right),\ |\lambda|\to\infty$, and analytic on $\partial D_+$. Thus, using Cauchy's theorem, the second integral in the RHS of \eqref{repre2} yields zero contribution.

Similarly, the sixth integral in the RHS of \eqref{repre2} vanishes. Indeed, since $\operatorname{Im}(\lambda) <0$  on $\partial D_-$, we get the asymptotic behaviour
\begin{equation*}
\begin{aligned}
\frac{1}{\Delta (\lambda)}  &\left[ (1- i \beta\nu)(1- i \alpha\lambda) \hat \theta (\lambda,t) -(1- i \beta\nu)(1- i \alpha\nu)\hat \theta (\nu,t)   \right] \\
&\sim  e^{i \nu L} \left[\frac{(1- i \beta\nu)(1- i \alpha\lambda)}{(1- i \beta\lambda)(1- i \alpha\nu)} \hat \theta (\lambda,t) -\frac{(1- i \beta\nu)}{(1- i \beta\lambda)}\hat \theta (\nu,t) \right],
\end{aligned} 
\end{equation*}
which is bounded and analytic on $\partial D_-$.

However, the two integrands might display singularities, namely the roots of the equation $\Delta (\lambda)=0$ which occur  in $D_\pm$ are poles in the corresponding domains of integrations. The position and the number of those poles depend on the values of $\alpha$ and $\beta$, but they are both always finite; see Appendix.
 In view of Remark \ref{ref-deform} there is always a favourable choice of $\partial D_\pm$ which allows the neglect of these poles.

We write \eqref{repre2} in compact form using the change of variables $\lambda \rightarrow \nu,$ for the integrals on $\partial D_-$,
which yields the result
\begin{equation*}
- \frac1{2\pi} \int_{\partial D_-} \frac{\phi (\lambda)}{\Delta (\lambda)}
  d\lambda = - \frac1{2\pi} \int_{\partial D_+} \frac{\phi (\nu)}{\Delta (\lambda)}
  d\lambda.
\end{equation*}

After some lengthy but straightforward calculations, we obtain
\begin{equation}\label{repre3}
\begin{aligned}
\theta (x,t) &= \frac1{2\pi} \int_{\R} e^{i \lambda x - \omega (\lambda)t} \hat \theta_0 (\lambda) d\lambda \\
&\phantom{=}- \frac{i}{\pi} e^{-\frac{K_0}{2 D_0}(L-x)} \int_{\partial D_{+}} \frac{e^{-\omega(\lambda) t}}{\Delta(\lambda)}  \left[ F_\alpha(\lambda, x)(1-i \beta \lambda) e^{i L\left(\lambda+i \frac{K_0}{2 D_0}\right)} \hat{\theta}_0(\lambda)\right.\\
&\phantom{=}\left.-F_{\beta}(\lambda, x-L)(1-i \alpha \nu) \hat{\theta}_0(\nu) +\left(K_0-i 2 \lambda D_0\right)\left( F_{\beta}(\lambda, x-L) \tilde{f} (\omega(\lambda), t) \right.\right.\\
&\phantom{=}\left.\left.-F_\alpha(\lambda, x) e^{-\frac{K_0}{2 D_0} L} \tilde{g}(\omega(\lambda),t)\right)\right] d \lambda,
\end{aligned}
\end{equation}
where
\begin{equation}\label{eq_f}
F_\gamma(\lambda, y)=\left(\gamma\frac{K_0}{2 D_0}-1\right) \sin \left(y\left(\lambda+i \frac{K_0}{2 D_0}\right)\right)-\gamma\left(\lambda+i \frac{K_0}{2D_0}\right) \cos \left(y\left(\lambda+i \frac{K_0}{2D_0}\right)\right),
\end{equation}
with 
$$
\Delta(\lambda)=e^{-i \lambda L}(1-i \alpha \lambda)(1-i \beta \nu )-e^{-i \nu L}(1-i \alpha \nu)(1-i \beta \lambda)
$$
and
$$
\omega(\lambda)=D_0 \lambda^2+i K_0 \lambda, \quad \nu=-\lambda-i \frac{K_0}{D_0} .
$$

From the above representation we can also deduce the solutions the following problems:
\begin{description}
\item[Robin--Dirichlet problem] We set $\beta=0$ in \eqref{repre3} which reduces to
\begin{equation}\label{repre4}
\begin{aligned}
\theta_{RD} (x,t) &= \frac1{2\pi} \int_{\R} e^{i \lambda x - \omega (\lambda)t} \hat \theta_0 (\lambda) d\lambda \\
&\phantom{=}- \frac{i}{\pi} e^{-\frac{K_0}{2 D_0}(L-x)} \int_{\partial D_{+}} \frac{e^{-\omega(\lambda) t}}{\Delta_\alpha(\lambda)}  \left\{ F_\alpha(\lambda, x) e^{i L\left(\lambda+i \frac{K_0}{2 D_0}\right)} \hat{\theta}_0(\lambda)\right.\\
&\phantom{=}\left.-(1-i \alpha \nu) \sin \left((L-x)\left(\lambda+i \frac{K_0}{2 D_0}\right)\right)\hat{\theta}_0(\nu) \right.\\
&\phantom{=}\left.+\left(K_0-i 2 \lambda D_0\right)\left[ \sin \left((L-x)\left(\lambda+i \frac{K_0}{2 D_0}\right)\right) \tilde{f} (\omega(\lambda), t) \right.\right.\\
&\phantom{=}\left.\left.-F_\alpha(\lambda, x) e^{-\frac{K_0}{2 D_0} L} \tilde{g}(\omega(\lambda),t)\right]\right\} d \lambda,
\end{aligned}
\end{equation}
where
\[
\Delta_\alpha (\lambda)=e^{-i \lambda L}(1-i \alpha \lambda)-e^{-i \nu L}(1-i \alpha \nu).
\]

\item[Dirichlet--Dirichlet problem] Setting $\alpha=\beta=0,$ the representation \eqref{repre3} coincides with \cite[Equation (25)]{arg24}.

\item[Neumann--Neumann problem] To derive the corresponding solution we set $\alpha=\beta$, and $f(t) = -\alpha f_N (t)$ and $g(t)=-\alpha g_N (t)$. Then, we take the limit $\alpha \rightarrow +\infty$, and the solution reads
\begin{equation*}
\begin{aligned}
\theta_{NN} (x,t) &= \frac1{2\pi} \int_{\R} e^{i \lambda x - \omega (\lambda)t} \hat \theta_0 (\lambda) d\lambda \\
&\phantom{=}+ \frac{1}{\pi} e^{-\frac{K_0}{2 D_0}(L-x)} \int_{\partial D_{+}} \frac{e^{-\omega(\lambda) t}}{\Delta_0 (\lambda)} \frac1{\lambda\nu}  \left\{ \lambda  G (\lambda, x)  e^{i L\left(\lambda+i \frac{K_0}{2 D_0}\right)} \hat{\theta}_0(\lambda)+\nu G(\lambda,x-L)\hat{\theta}_0(\nu) \right.\\
&\phantom{=}\left.+i \left(K_0-i 2 \lambda D_0\right)\left[ G(\lambda,x-L) \tilde{f}_N (\omega(\lambda), t) -G(\lambda, x) e^{-\frac{K_0}{2 D_0} L} \tilde{g}_N (\omega(\lambda),t)\right]\right\} d \lambda,
\end{aligned}
\end{equation*}
where
\[
G(\lambda,y) = \lim_{\alpha\rightarrow +\infty} \frac{F_\alpha (\lambda,y)}{a} =  \frac{K_0}{2 D_0} \sin \left(y\left(\lambda+i \frac{K_0}{2 D_0}\right)\right)-\left(\lambda+i \frac{K_0}{2D_0}\right) \cos \left(y\left(\lambda+i \frac{K_0}{2D_0}\right)\right)
\]
and
\[
\Delta_0 (\lambda)=e^{-i \lambda L}-e^{-i \nu L}.
\]

\end{description}

\subsection{Numerical evaluation}\label{sec_num_direct}

In what follows we present two examples of numerical evaluation of the solution \eqref{repre3}, for two different physical setups; our solutions match with the ones existing in the literature. We find worth mentioning that, similarly to the discussion in \cite{arg24}, we derive and compute the general solution through a unified methodology,  whereas the existing solutions were obtained through different approaches, which sometimes only yield approximate results. Furthermore, the (integral) form of \eqref{repre3} is favourable for numerical evaluation compared to existing series representations.

\subsubsection{Robin--Dirichlet problem}

In \cite{Bra73} Braester considers the problem of vertical infiltration at constant flux $q$ into a bounded medium, with length $L,$ with a constant water table. The governing IBVP reads:
\begin{equation*}
\begin{aligned}
\frac{\partial \theta}{\partial t} + K_0 \frac{\partial \theta}{\partial x} &= D_0 \frac{\partial^2 \theta}{\partial x^2},  \quad && 0<x <L, \, t>0, \\
\theta (x,0) &= \theta_0, \quad && 0 <x <L, \\ 
D_0 (\theta (0,t)-\theta_0) - D_0\frac{\partial \theta}{\partial x} (0,t) &= q,  \quad &&t>0,\\
\theta (L,t)  &= \theta_s,  \quad &&t>0, 
\end{aligned}
\end{equation*}
where $\theta_0$ and $\theta_s$ are the initial and saturation values of the water content, respectively.

We define $u(x,t) = \theta (x,t) - \theta_0$ that solves
\begin{equation*}
\begin{aligned}
\frac{\partial u}{\partial t} + K_0 \frac{\partial u}{\partial x} &= D_0 \frac{\partial^2 u}{\partial x^2},  \quad && 0<x <L, \, t>0,  \\
u (x,0) &= 0, \quad && 0 <x <L, \\ 
u (0,t) - \frac{\partial u}{\partial x} (0,t) &= \frac{q}{D_0},  \quad &&t>0, \\
u (L,t)  &= \theta_s - \theta_0,  \quad &&t>0. 
\end{aligned}
\end{equation*}

The general solution is given by \eqref{repre4}. Given the homogeneous initial condition, the constant boundary functions and that $\alpha=1,$ we get the simplified form
\begin{equation}\label{ex1sol}
\begin{aligned}
\theta (x,t) &= \theta_0 - \frac{i}{\pi} e^{-\frac{K_0}{2 D_0}(L-x)} \int_{\partial D_{+}} \frac{1-e^{-\omega(\lambda) t}}{\omega(\lambda)\Delta_1(\lambda) }  \left(K_0-i 2 \lambda D_0\right)\left[\frac{q}{D_0} \sin \left((L-x)\left(\lambda+i \frac{K_0}{2 D_0}\right)\right)  \right.\\
&\phantom{=}\left.- (\theta_s - \theta_0) F_1(\lambda, x) e^{-\frac{K_0}{2 D_0} L} \right]d \lambda,
\end{aligned}
\end{equation}
where
\[
\Delta_1 (\lambda)=e^{-i \lambda L}(1-i  \lambda)-e^{-i \nu L}(1-i  \nu)
\]
and
\[
F_1 (\lambda, y)=\left(\frac{K_0}{2 D_0}-1\right) \sin \left(y\left(\lambda+i \frac{K_0}{2 D_0}\right)\right)-\left(\lambda+i \frac{K_0}{2D_0}\right) \cos \left(y\left(\lambda+i \frac{K_0}{2D_0}\right)\right).
\]
Here we have used that $\hat u_0 (\lambda) = \hat u_0 (\nu) =0$ because of the homogeneous initial condition and that (see \eqref{t_transform})
\[
\tilde{f} (\omega(\lambda), t) = \frac{q}{D_0} \frac{e^{\omega(\lambda) t}-1}{\omega(\lambda)}, \quad \tilde{g} (\omega(\lambda), t) = (\theta_s - \theta_0) \frac{e^{\omega(\lambda) t}-1}{\omega(\lambda)}.
\]
Regarding the evaluation of the integral in \eqref{ex1sol} we observe that the roots of the equation $\omega(\lambda)\Delta_1(\lambda) =0$ are located in the lower half complex plane $\Im (\lambda) \leq 0$ and thus we can compute the integral by simply deforming $\partial D_{+}$ so that it lies in the upper half complex plane. For a thorough discussion see \cite{arg24}.

\begin{figure}
\centering
\begin{subfigure}{.5\textwidth}
  \centering
\includegraphics[width=0.8\textwidth]{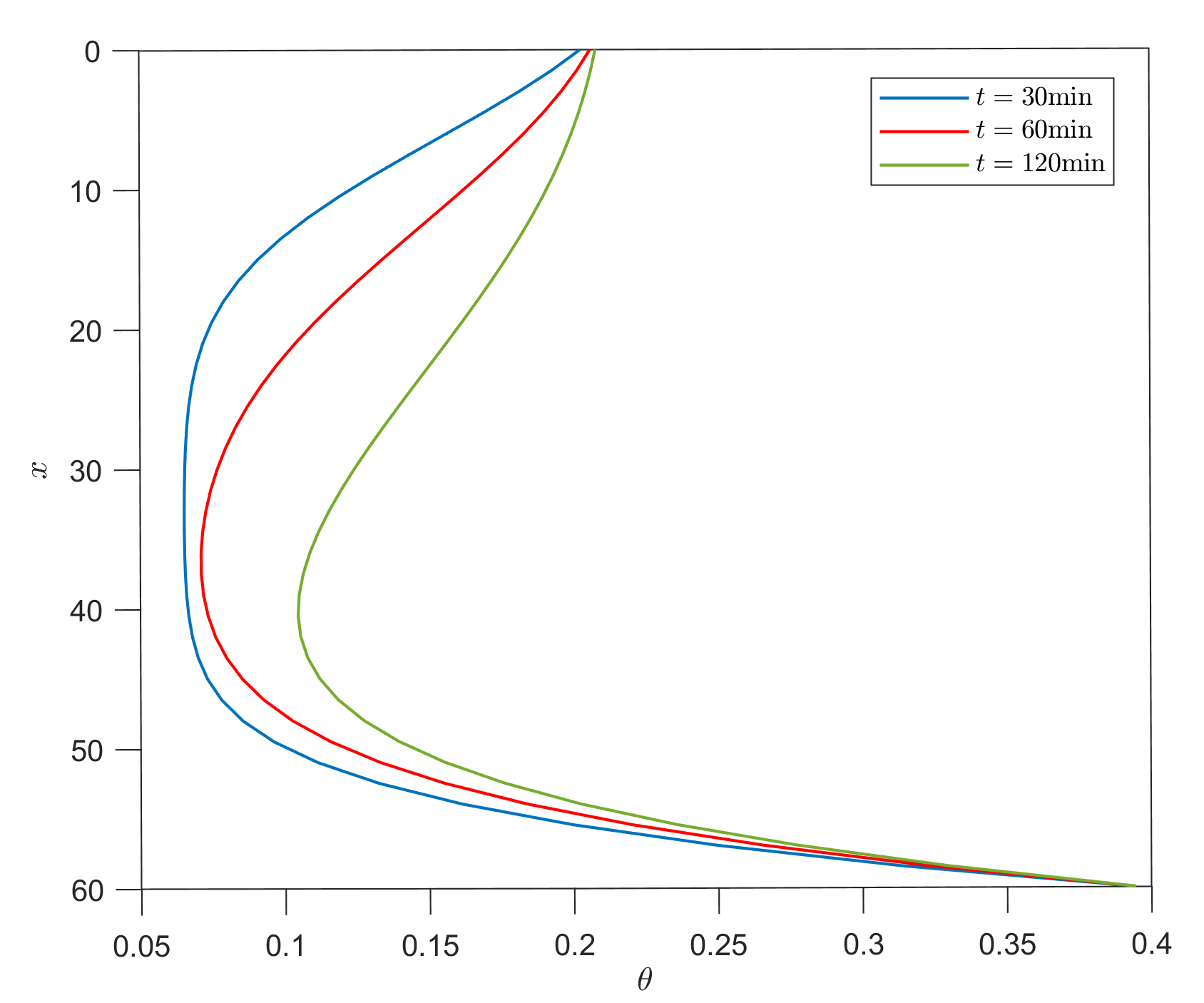}
  \label{fig:sub1}
\end{subfigure}%
\begin{subfigure}{.5\textwidth}
  \centering
\includegraphics[width=0.9\textwidth]{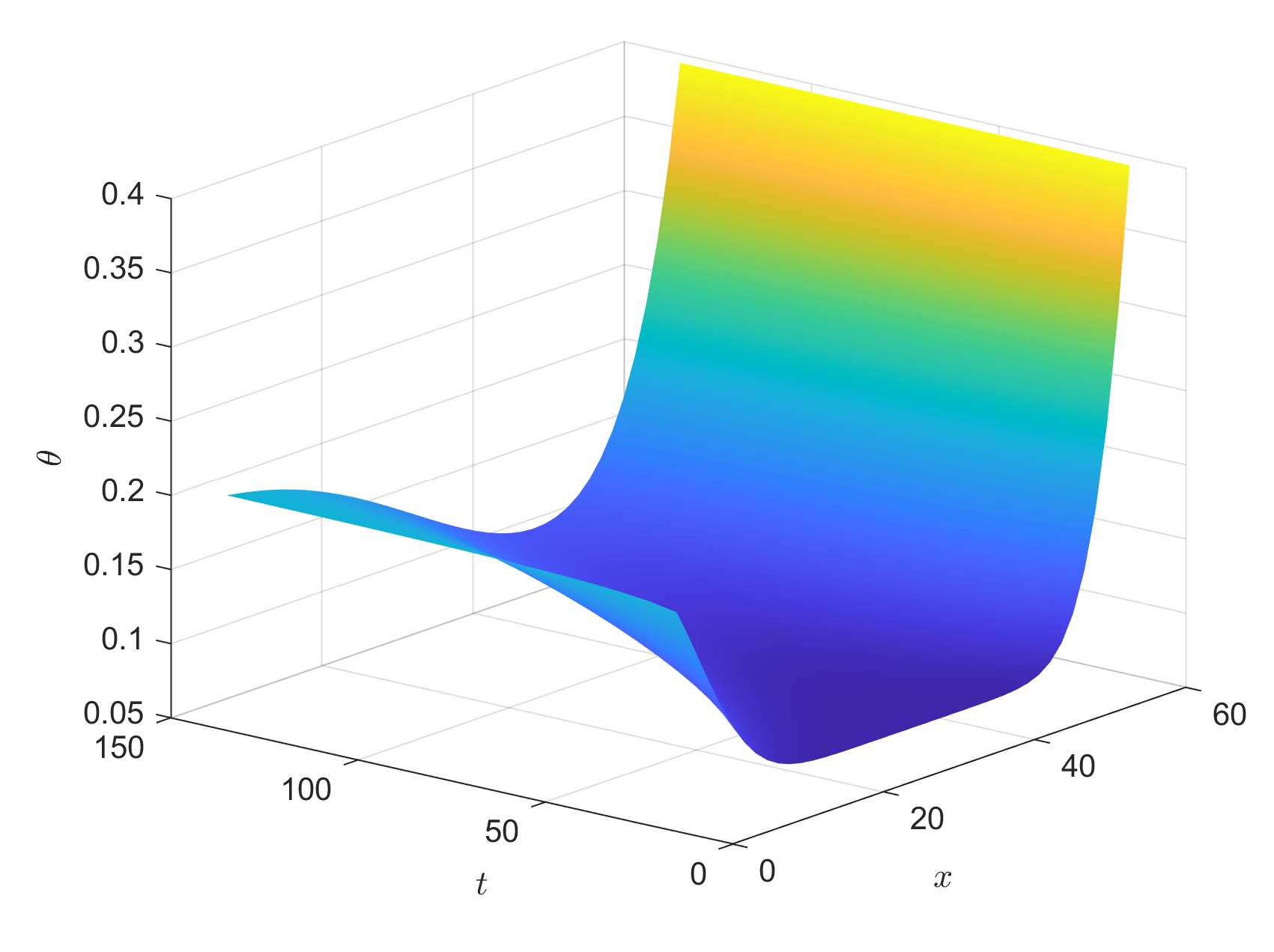}
  \label{fig:sub2}
\end{subfigure}
\caption{The water content $\theta$ given by \eqref{ex1sol} at fixed time steps (left) and for $x \in [0,\,L]$ and $t\in[0,\,135]$min (right).}
\label{fig12}
\end{figure}

As soil we consider the `Rehovot Sand' with parameters $D_0 = 0.208\times 10^{-1} \mbox{cm}^2/\mbox{sec},$ $K_0 = 0.144 D_0,$ and $\theta_0 = 65\times 10^{-3}$ and $\theta_s= 397\times 10^{-3}$  \cite{Bra73}. The constant flux at the surface $x=0$ is given by $q=0.3 \times 10^{-3}\mbox{cm}/\mbox{sec}$ and the total length is $L=60 \mbox{cm}.$ In the left picture of \autoref{fig12} we present the solution derived by \eqref{ex1sol} at the time steps $t = 30, \, 60$ and $120$min. The solution follows the expected behaviour and it is qualitatively identical with the one presented in \cite[Fig. 4]{Bra73}. In the right-hand side of \autoref{fig12} the solution is plotted for $x \in [0,\,L]$ and $t\in[0,\,135]$min.

\subsubsection{Robin--Robin problem}

The study of the initial problem \eqref{bvp} helps us to consider also quasi-linear soils where the hydraulic conductivity $K$ and the diffusivity $D$ are arbitrary functions of the water content but connected through $\tfrac{d K}{d \theta} = c D,$ for $c>0$ \cite{Philip91}. Thus, the problem of constant rate rainfall infiltration into bounded shallow profiles can be modelled with the following IBVP for the conductivity: 
\begin{equation*}
\begin{aligned}
\frac{\partial K}{\partial t} + \frac{\partial K}{\partial x} &= \frac12 \frac{\partial^2 K}{\partial x^2},  \quad && 0<x <L, \, t>0,  \\
K (x,0) &= 0, \quad && 0 <x <L, \\ 
K (0,t) -\frac12 \frac{\partial K}{\partial x} (0,t) &= R,  \quad &&t>0, \\
K (L,t) -\frac12 \frac{\partial K}{\partial x} (L,t) &= 0,  \quad &&t>0, \
\end{aligned}
\end{equation*}
where $R>0$ is the flow velocity. Note that in this formulation $t$ and $x$ are not the physical time and spatial coordinate, respectively, but normalized quantities, see \cite{philip1989} for more details.

Considering \eqref{repre3} for $K_0 = 1, \, D_0 = \tfrac{1}2$ and $\alpha=\beta=\tfrac{1}2,$ and the above initial and boundary conditions, we obtain
\begin{equation}\label{sol_robin}
K (x,t) = - R\frac{i}{\pi} e^{x-L} \int_{\partial D_{+}} \frac{1-e^{-\omega(\lambda) t}}{\omega(\lambda)\Delta(\lambda)}  \left(1-i  \lambda \right) F_{1/2}(\lambda, x-L)  d \lambda,
\end{equation}
where $F$ is given by \eqref{eq_f} for $K_0 = 1, \, D_0 = \tfrac{1}2$ and 
\[
\Delta (\lambda)= \left(1-  \frac{i \lambda}2\right) \left(1-\frac{i \nu}2\right) \left(e^{-i \lambda L}-e^{-i \nu L} \right).
\]

\begin{figure}
\centering
\begin{subfigure}{.5\textwidth}
  \centering
\includegraphics[width=0.8\textwidth]{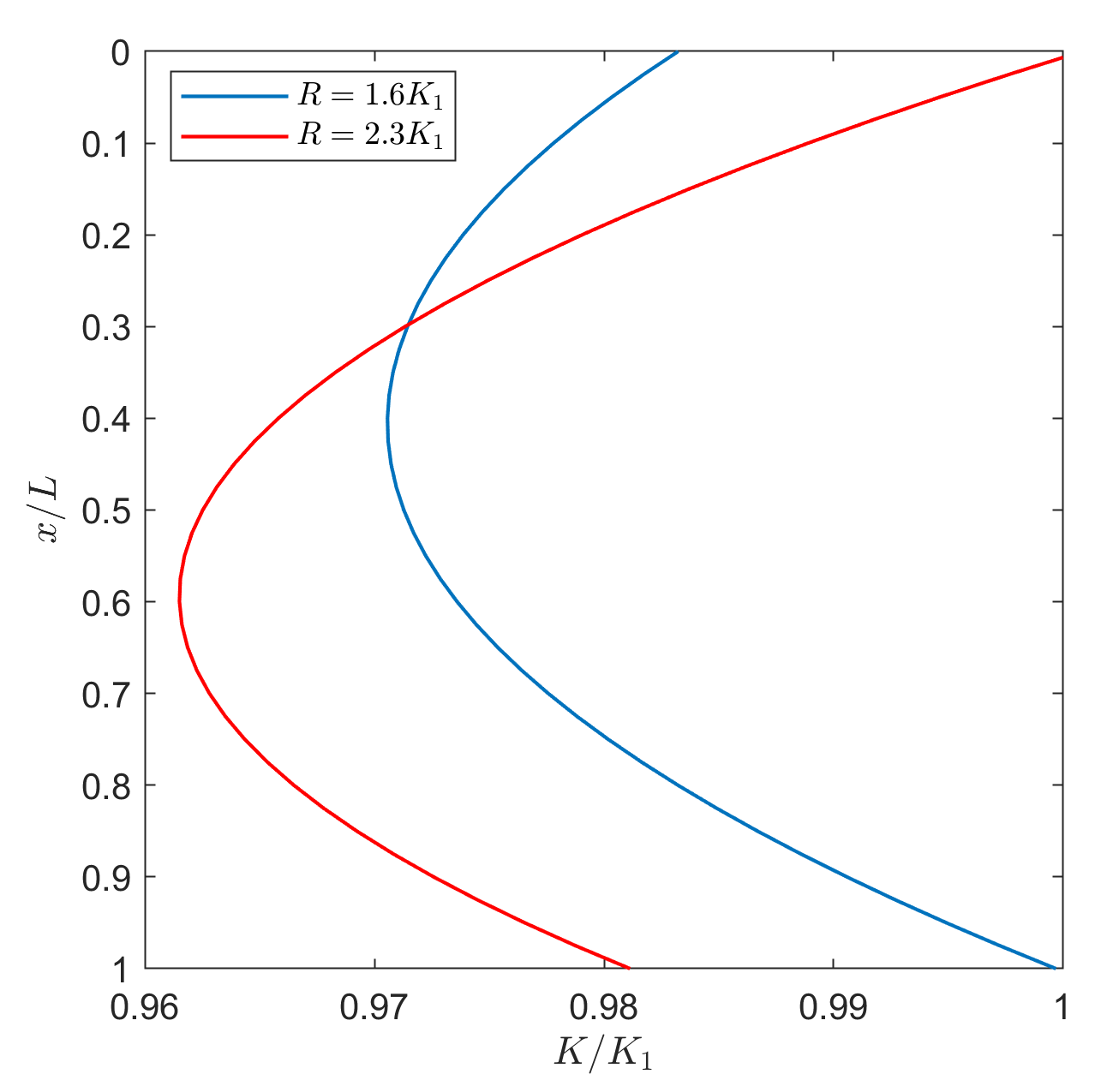}
\end{subfigure}%
\begin{subfigure}{.5\textwidth}
  \centering
\includegraphics[width=1\textwidth]{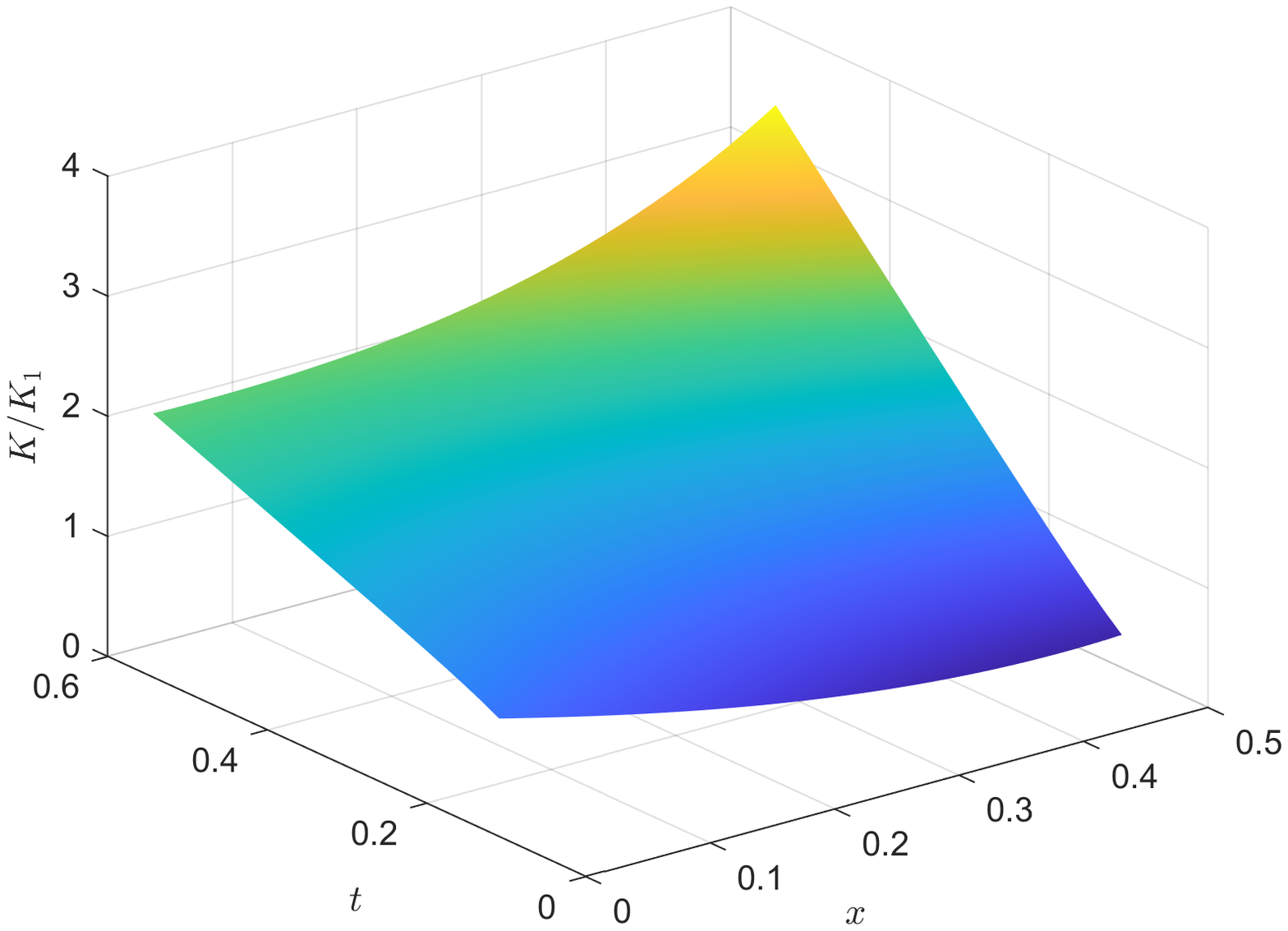}
\end{subfigure}
\caption{The conductivity $\tfrac{K}{K_1}$ given by \eqref{sol_robin} for different flow velocities at fixed time $t=\tfrac{L K_1}R$ (left) and for $x \in [0,\,L]$ and $t\in[0.5,\,2.5] \tfrac{L K_1}R$ (right).}
\label{fig34}
\end{figure}

To reproduce the results from \cite[Fig. 4]{Philip91} we consider two different instances during infiltration for $R=1.6 K_1$ (base saturation) and for $R=2.3 K_1$ (surface ponding) where $K_1$ is the saturated hydraulic conductivity. The profiles of $K$ for $L=0.05$ at $t=\tfrac{L K_1}R$ are shown in \autoref{fig34}. Again we get a perfect match. In the right picture of \autoref{fig34} we present the solution for $R=2.3 K_1$ and $L=0.5.$ The time interval is $t\in \left[\tfrac{L K_1}{2R},\,\tfrac{5L K_1}{2R}\right].$

\section{Inverse/Control Problem}\label{sec_inv}

In this section we examine the inverse problem of finding the boundary condition at $x=L$ that results in a steady state (i.e. vanishing) solution at a final time $t=T.$ In particular, we consider the IBVP:
\begin{equation*}
\begin{aligned}
\frac{\partial \theta}{\partial t} + K_0 \frac{\partial \theta}{\partial x} &= D_0 \frac{\partial^2 \theta}{\partial x^2},  \quad && 0<x <L, \, t>0,  \\
\theta (x,0) &= \theta_0 (x), \quad && 0 <x <L, \\ 
\theta (0,t) -\alpha \frac{\partial \theta}{\partial x} (0,t) &= 0,  \quad &&t>0, \\
\theta (L,t) &= v(t),  \quad &&t>0, 
\end{aligned}
\end{equation*}
and we examine the following boundary control problem: Given the initial condition $\theta_0,$ find the boundary control $v(t)$ such that $\theta(x,T)=0,$ for $0<x<L.$ In \cite{hwang23} a similar problem was considered for $a=0$ and control at the boundary $x=0$, based on the method introduced in \cite{KOD23} where the null controllability was examined for the heat equation ($K_0 =0$) and Neumann-Neumann boundary conditions.

If $\theta (x,T)=0$ holds true, then from \eqref{repre4} we get
\begin{equation}\label{control_eq}
- \frac{i}{\pi} e^{-\frac{K_0}{2 D_0}(2L-x)} \int_{\partial D_{+}} \frac{e^{-\omega(\lambda) T}}{\Delta_\alpha(\lambda)}  \left(K_0-i 2 \lambda D_0\right) F_\alpha(\lambda, x)  \tilde{v}(\omega(\lambda),T) d \lambda = B(x,T),
\end{equation}
where
\begin{equation*}
\begin{aligned}
B (x,t) &= \frac1{2\pi} \int_{\R} e^{i \lambda x - \omega (\lambda)t} \hat \theta_0 (\lambda) d\lambda \\
&\phantom{=}- \frac{i}{\pi} e^{-\frac{K_0}{2 D_0}(L-x)} \int_{\partial D_{+}} \frac{e^{-\omega(\lambda) t}}{\Delta_\alpha(\lambda)}  \left[ F_\alpha(\lambda, x) e^{i L\left(\lambda+i \frac{K_0}{2 D_0}\right)} \hat{\theta}_0(\lambda)\right.\\
&\phantom{=}\left.-(1-i \alpha \nu) \sin \left((L-x)\left(\lambda+i \frac{K_0}{2 D_0}\right)\right)\hat{\theta}_0(\nu) \right] d \lambda.
\end{aligned}
\end{equation*}

Following \cite{KOD23} we assume a series representation of $v(t)$ for a sine-Fourier basis supported in $(\tau, T),$ for $\tau>0.$ The truncated series approximation reads
\begin{equation}\label{eq_control}
v(t)=\sum_{n=1}^{N+1} c_n \phi_n(t),
\end{equation}
where
$$
\phi_n(t)= \begin{cases} 0, &  0 \leq t \leq \tau, \\
 \sin \left(n \pi \dfrac{t-\tau}{T-\tau}\right), &  \tau \leq t \leq T, \\ 
 0, & t>T.\end{cases}
$$
Then, using \eqref{t_transform} we get
\begin{equation}\label{eq_vhat}
\tilde{v}(\omega(\lambda), t)=\sum_{n=1}^{N-1} c_n \varphi_n(\lambda, t),
\end{equation}
where
$$
\begin{aligned}
\varphi_n(\lambda, t) &:= \int_0^t e^{\omega(\lambda) s } \phi_n(s) d s
\\
&= \frac{\tau-T}{\left(\tau-T\right)^2 \omega(\lambda)^2+n^2 \pi^2} \left[-n \pi e^{\omega(\lambda) \tau}+e^{\omega(\lambda) t}\left(n \pi \cos \left(\frac{n \pi\left(\tau-t\right)}{\tau-T}\right)\right.\right. \\
&\left.\left.+\left(\tau-T\right) \omega(\lambda) \sin \left(\frac{n \pi\left(\tau-t\right)}{\tau-T}\right)\right)\right], \quad \tau \leq t \leq T
\end{aligned}
$$
and $\varphi_n(\lambda, t)=0$ for $0<t<\tau,$ or $t>T.$

Substituting \eqref{eq_vhat} for $t=T$ in \eqref{control_eq} we derive the linear equation
\begin{equation}\label{eq_system}
\sum_{n=1}^{N+1} c_n A_n(x, T)= B(x, T), \quad x\in(0,L),
\end{equation}
where
\[
A_n (x,t) =  - \frac{i}{\pi} e^{-\frac{K_0}{2 D_0}(2L-x)} \int_{\partial D_{+}} \frac{e^{-\omega(\lambda) t}}{\Delta_\alpha(\lambda)}  \left(K_0-i 2 \lambda D_0\right) F_\alpha(\lambda, x)  \varphi_n(\lambda, t) d \lambda.
\]

To solve \eqref{eq_system} we discretize the domain $(0, L)$ considering equidistant grid points $x_k,\, k = 0,\ldots,N.$ Then we obtain $N+1$ equations of the form
\begin{equation*}
\sum_{n=1}^{N+1} c_n A_{nk}( T)= B_k ( T), \quad k=0,\ldots,N,
\end{equation*}
where $A_{nk}( T) = A_{n}(x_k, T)$ and $B_{k}( T)= B(x_k, T).$ The above well-determined linear system of $N+1$ equations for the $N+1$ unknown coefficients of $v$ reads in compact form
\begin{equation}\label{eq_system3}
\bb A \bb c = \bb b,
\end{equation}
where $\bb A$ is $(N+1)\times (N+1)$ matrix and $\bb b,\, \bb c = [c_1, \, c_2, \ldots, c_{N+1} ]^\top$ are 
$(N+1)\times 1$ vectors. All the involved integrals are computed numerically and we derive $\bb c= \bb A^{-1} \bb b$ to be substituted in \eqref{eq_control} to obtain the boundary control.

\subsection{Numerical examples}

In the following examples, we set the control to act immediately, namely $\tau=0$ and we choose the points of discretisation $x_k = \tfrac{(k+1)L}{N+2}, \, k = 0,...,N$. We consider four examples in order to illustrate some features of the algorithm:
\begin{itemize}
\item[1.] First we consider a typical problem, with discontinuous initial condition, which illustrates the exponential decay of the error (deviation from the desired state), in terms of the discretisation points, whereas the order of magnitude of the control itself remains the same.
\item[2.] In the case of a small diffusivity coefficient $D_0>0$ the situation is qualitatively identical, but the control inevitably grows.  We revisit this problem, by employing a classical regularisation scheme, which mitigates this effect by relaxing the constraint on the deviation between the desired and the achieved final state.
\item[3.] We consider, as the special case $K_0=\alpha=0$, the Dirichlet problem of the heat equation. We show that this simple problem can be null-controlled with very little effort, namely with only 3 points of discretisation we obtain a control which is 10 times smaller than the initial condition (for an exponentially smaller error).
\item[4.] Using the setup of example 3, we study the effect when very little time is allowed for the action of the control, namely when $T>0$ is small. The regularisation process acts similarly to example 2.
\end{itemize}

These results indicate that large diffusivity $D_0>0$ would act in favour of null-controlling the IBVP. In other words, smaller diffusivity requires much  larger controlling effort for given time $T>0$. The same effect is observed in the case that diffusivity remains the same, but $T>0$ becomes very small. This observation indicates an obvious correspondence between these two setups. In both cases, the regularisation scheme allows for much smaller controls.

\textbf{Example 1:} In the first example we consider a piecewise constant initial condition of the form
\begin{equation}\label{initial_piecewise}
\theta_0 (x) = \left.
  \begin{cases}
    1, & \text{for } x \in \left(0,\tfrac{L}2\right), \\
    0, & \text{for } x \in \left(\tfrac{L}2,L\right)
  \end{cases}
  \right.
\end{equation}
and considering \eqref{fourier} we get
\[
\hat  \theta_0 (\lambda) = \frac{-i}{\lambda} \left( 1- e^{-i \lambda\tfrac{L}{2}} \right).
\]

\begin{figure}[h!]
\begin{center}
\includegraphics[width=0.8\textwidth]{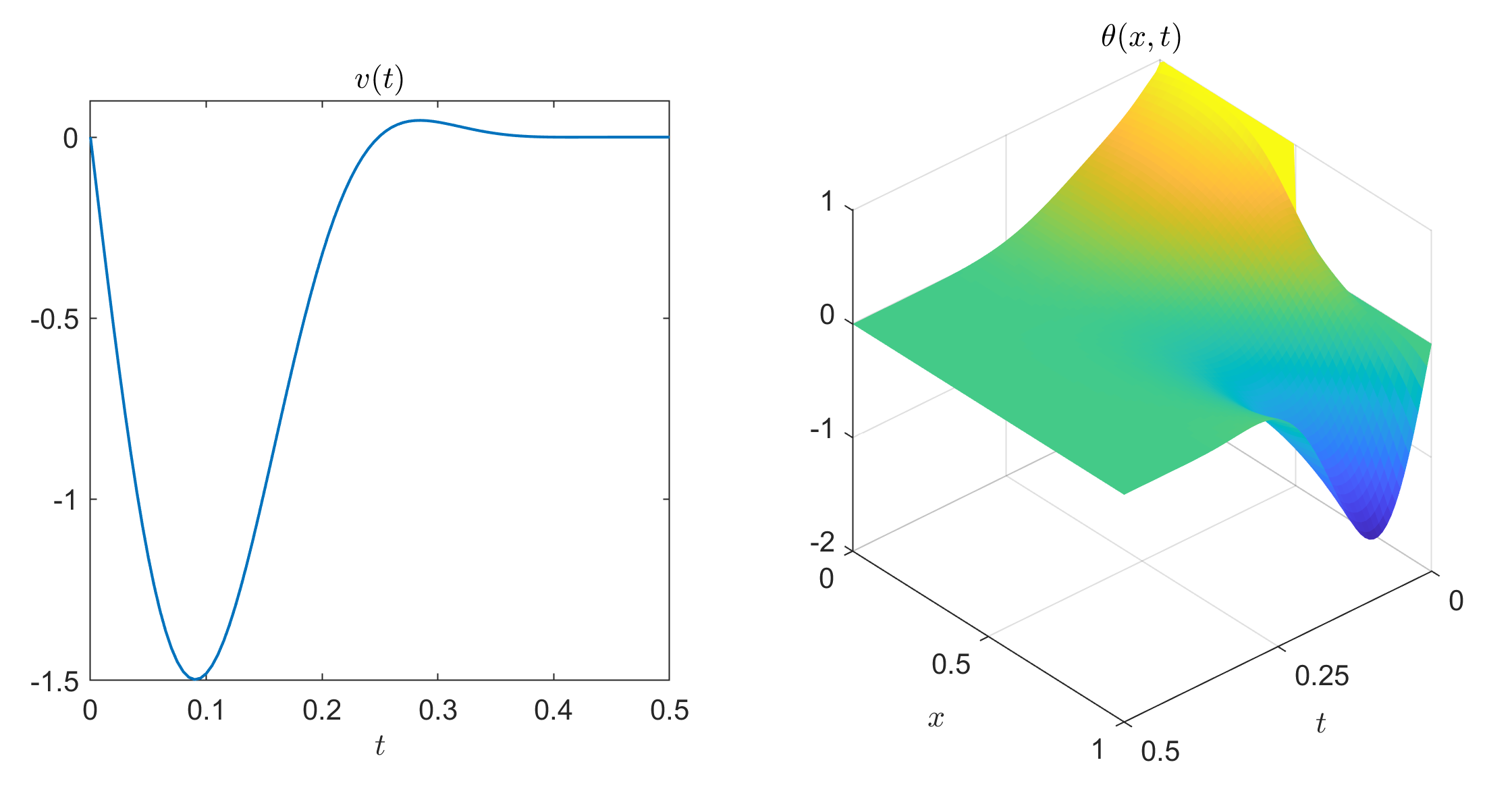}
\end{center}\caption{The control $v(t)$ (left) and the solution $\theta(x,t)$ (right) for $N=4$ of the first example.}\label{fig_ex1_3d}
\end{figure}

We set $D_0 =1$ and $K_0 = 1/2,$ for $L=1$ and $a=1.$ By varying the final time $T$ we want to examine the convergence of the solution with respect to $N.$ In \autoref{table1} we present the $L^2$ norm of the control $v(t)$ and in \autoref{table2} the $L^2$ norm of the vanishing solution $\theta (x,T)$, namely the deviation from the desired state. For all values of $T$ we obsevre exponential decay of $\theta (x,T)$ with $N$, whereas the control $v(t)$ retains the same order of magnitude.
 The control and the solution for $N=4$ and $T=0.5$ are presented in \autoref{fig_ex1_3d}.

\begin{table}[h!]
\begin{center}
 \begin{tabular}{ | c  | c  | c  | c  | } 
 \hline
 $N \mathbin{\bigg\backslash} T$ & $1/2$ & $1$ & $2$  
\\ \hhline{====}
  2    & 0.346233  & 0.108330  & 0.012763
 \\
   4    &  0.48019 & 0.170745 & 0.030665
  \\ 
   6    &  0.595987 & 0.223610 & 0.050358
    \\ \hline
 \end{tabular}
\caption{The $L^2$ norm of the control $v(t)$ of the first example for $D_0=1$ and initial data \eqref{initial_piecewise}.}\label{table1}
\end{center}
\end{table}

\begin{table}[h!]
\begin{center}
 \begin{tabular}{ | c  | c  | c  | c  | } 
 \hline
 $N \mathbin{\bigg\backslash} T$ & $1/2$ & $1$ & $2$  
\\ \hhline{====}
  2    & 5.816169E-5 & 4.749178E-7 & 1.290414E-9
 \\
   4    & 1.621817E-8 & 1.012742E-11 & 2.673885E-15
  \\ 
   6    & 1.340275E-12 & 5.966348E-17 & 1.262254E-21
    \\ \hline
 \end{tabular}
\caption{The $L^2$ norm of the solution $\theta(x,T)$ of the first example for $D_0 =1$ and initial data \eqref{initial_piecewise}.}\label{table2}
\end{center}
\end{table}

\textbf{Example 2:} In this example,  the initial condition is given by
\begin{equation}\label{initial_cos}
\theta_0 (x) = \cos \left(\frac{\pi}{2L}x \right), \quad x\in (0,L),
\end{equation}
resulting in
\[
\hat  \theta_0 (\lambda) = -\frac{2L}{4 \lambda^2 L^2 - \pi^2} \left( 2 i \lambda L+ \pi e^{-i \lambda L} \right).
\]

We keep all parameters the same except for the diffusivity which is now  $D_0 = 0.1.$ In \autoref{table3} and \autoref{table4} we present the values of the control and the solution at the final time, respectively. The cases $N=12,\,14,\,16$ and $T=0.5, \, 1,\,2$ are considered.

\begin{table}[h!]
\begin{center}
 \begin{tabular}{ | c  | c  | c  | c  | } 
 \hline
 $N \mathbin{\bigg\backslash} T$ & $1/2$ & $1$ & $2$  
\\ \hhline{====}
  12    & 198369 & 2746.9 & 195.0
 \\
   14    &   384725 & 4020.3 & 242.0
  \\ 
   16    &  702459 & 5702.6 & 295.5
    \\ \hline
 \end{tabular}
\caption{The $L^2$ norm of the control $v(t)$ of the second example for $D_0 = 0.1$ and initial data \eqref{initial_cos}.}\label{table3}
\end{center}
\end{table}

\begin{table}[h!]
\begin{center}
 \begin{tabular}{ | c  | c  | c  | c  | } 
 \hline
 $N \mathbin{\bigg\backslash} T$ & $1/2$ & $1$ & $2$  
\\ \hhline{====}
  12    & 1.372320E-2  & 3.754687E-9 & 3.245673E-16
 \\
   14    &  3.842986E-4 & 8.502099E-12 & 5.342716E-20
  \\ 
   16    &  6.815284E-6 & 1.195998E-14 & 5.400493E-24
    \\ \hline
 \end{tabular}
\caption{The $L^2$ norm of the solution $\theta(x,T)$ of the second example for $D_0 = 0.1$ and initial data \eqref{initial_cos}.}\label{table4}
\end{center}
\end{table}

\begin{figure}[h!]
\begin{center}
\includegraphics[width=0.8\textwidth]{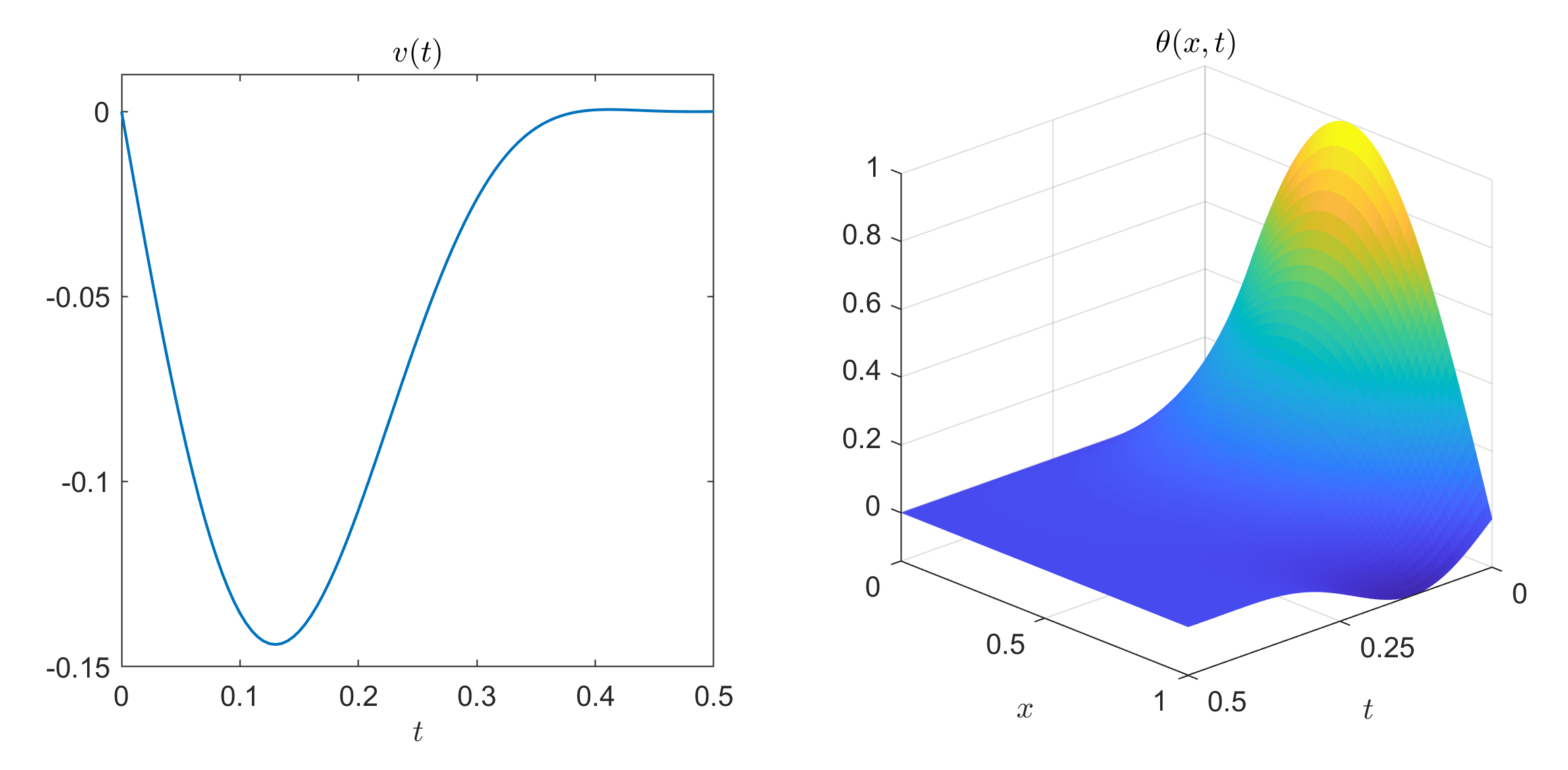}
\end{center}\caption{The control $v(t)$ (left) and the solution $\theta(x,t)$ (right) for $N=2$ of the third example.}\label{fig_ex3_3d}
\end{figure}

\textbf{Example 3:}
We consider the Dirichlet-Dirichlet problem ($\alpha=0$) with solution given by \cite[Equation (25)]{arg24}. Let now
\begin{equation}\label{initial_sin}
\theta_0 (x) = \sin \left(\frac{\pi}{L}x \right), \quad x\in (0,L),
\end{equation}
to obtain
\[
\hat  \theta_0 (\lambda) = -\frac{L \pi}{\lambda^2 L^2 - \pi^2} \left(1+  e^{-i \lambda L} \right).
\]
In this example we set $K_0 =0,$ meaning we consider the heat equation, and  $D_0=1.$  The control and the solution for $N=2$ and $T=1/2$ are presented in \autoref{fig_ex3_3d} and their norms for varying $N$ and $T$ in \autoref{table5} and \autoref{table6}, respectively.

\begin{table}[h!]
\begin{center}
 \begin{tabular}{ | c  | c  | c  | c  | } 
 \hline
 $N \mathbin{\bigg\backslash} T$ & $1/2$ & $1$ & $2$  
\\ \hhline{====}
  2    & 5.527160E-2 & 3.094840E-3  & 3.353640E-6
 \\
   4    &   9.95742E-2  & 1.004530E-2  & 4.296490E-5
  \\ 
   6    &  1.409660E-1 & 1.963610E-2 & 2.074790E-4
    \\ \hline
 \end{tabular}
\caption{The $L^2$ norm of the control $v(t)$ of the third example for $K_0 = 0$ and initial data \eqref{initial_sin}.}\label{table5}
\end{center}
\end{table}

\begin{table}[h!]
\begin{center}
 \begin{tabular}{ | c  | c  | c  | c  | } 
 \hline
 $N \mathbin{\bigg\backslash} T$ & $1/2$ & $1$ & $2$  
\\ \hhline{====}
  2    & 1.752254E-6  & 2.309529E-9 & 5.618245E-14
 \\
  4    &   3.517729E-10  & 5.396213E-14 & 3.264470E-19
  \\ 
  6    &  2.144497E-14  & 2.944358E-19 & 2.774625E-25
    \\ \hline
 \end{tabular}
\caption{The $L^2$ norm of the solution $\theta(x,T)$ of the third example for $K_0 = 0$ and initial data \eqref{initial_sin}.}\label{table6}
\end{center}
\end{table}

\subsubsection{Regularized controls}

As it is already mentioned (very) small diffusivity, $D_0>0$, and/or action time, $T>0$, could result in a large control in order to achieve high accuracy in the solution of the inverse problem (see e.g.  \autoref{table3}).
To handle this problem we intend to determine the control, namely the vector of coefficients $\bb c$ by  \eqref{eq_system3}, considering a regularization scheme rather than solve it directly. The main idea is to sacrifice/relax the error of the solution in order to achieve a smaller control. We formulate a constrained optimization problem of the form
\[
\min \, \lVert \bb c\rVert_2  \quad \mbox{subject to} \quad \lVert \bb A  \bb c - \bb b\rVert_2 \leq \delta,
\]
for some reasonably small $\delta>0$. The idea is that $\delta$ is large enough to allow the relaxing of the restriction of the `exact' problem, namely it would be larger than the error induced by the  numerical solution of \eqref{eq_system3}, which is given as the $L^2$ norm of the `exact' solution $\theta (x,T).$
 In the following examples $\delta$ is of the order of $10^{-3}.$ In other words, the parameter $\delta$ plays the role of the regularization parameter and the solution of the above least squares problem is similar to the regularized solution from classical Tikhonov regularization for appropriately chosen values of the regularization parameter \cite{hansen2010}. This part of the work is far from being an exhaustive investigation of regularisation methods. It has to be seen as an initial step towards regularized controls; more sophisticated schemes will be examined in future work. However, valuable insights can be drawn despite the classical regularization tools used. 

\begin{figure}
\begin{center}
\includegraphics[width=0.9\textwidth]{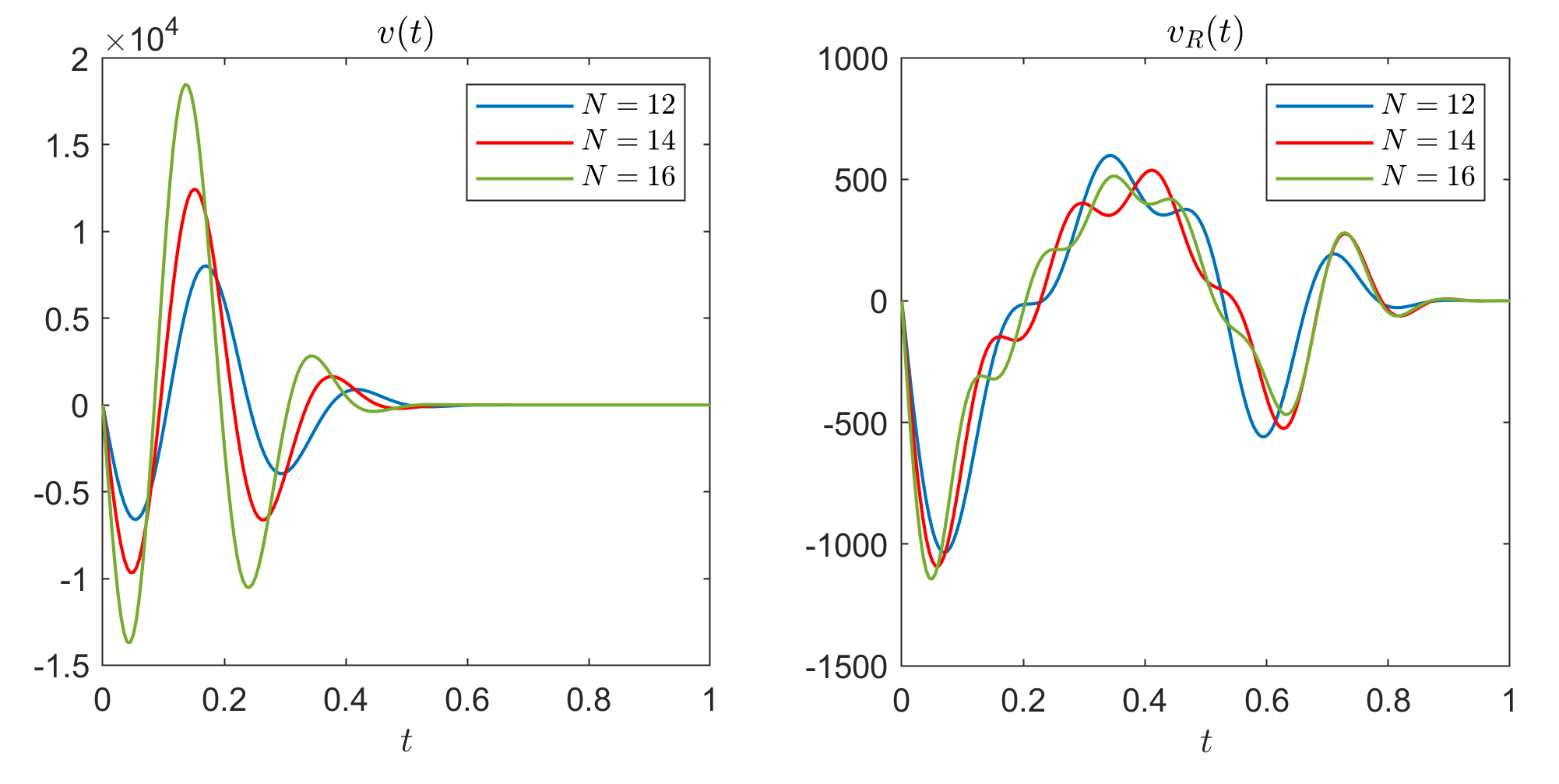}
\end{center}\caption{The exact (left) and the regularized (right) controls of the second example for different values of $N$.  }\label{fig5}
\end{figure}

\textbf{Example 2 (regularized)} We  examine the effect of the aforementioned regularization scheme, by revisiting example 2 for $T=1$ and $N=12,\,14$ and $16$. We choose $\delta= O\left( 10^{-3} \right)$ so that the error of the solution is of the same order for all values of $N$. 
The controls are depicted in \autoref{fig5}; the magnitude of the regularized control $v_R$, as well as the norm of the corresponding error/solution $\theta_R$ are presented in \autoref{table7}.

\begin{table}[ht]
\begin{center}
 \begin{tabular}{ | c  | c  | c  |c|| c | c|} 
 \hline
 $N$ & $\lVert v(t)\rVert_2$ &  $\lVert v_R(t)\rVert_2$ & Rel. change & $\lVert\theta(x,T)\rVert_2$  & $\lVert\theta_R(x,T)\rVert_2$  
\\ \hhline{======}
12 & 2746.9 &  388.5  & $-85.9\%$ &  3.754687E-9 & 5.237347E-3  
 \\
   14 & 4020.3 & 373.6  & $-90.7\%$ &  8.502099E-12 &   9.039701E-3  
  \\ 
   16 & 5702.6 &  370.5  & $-93.5\%$ & 1.195998E-14 &  6.922691E-3   
    \\ \hline
 \end{tabular}
\caption{The second example for $D_0 =0.1$ and $T=1$: The first two columns are the $L^2$ norms of the exact and regularized control. The last two columns are $L^2$ norms of the associated errors, namely the solutions at the final time.}\label{table7}
\end{center}
\end{table}

\begin{figure}
\begin{center}
\includegraphics[width=0.9\textwidth]{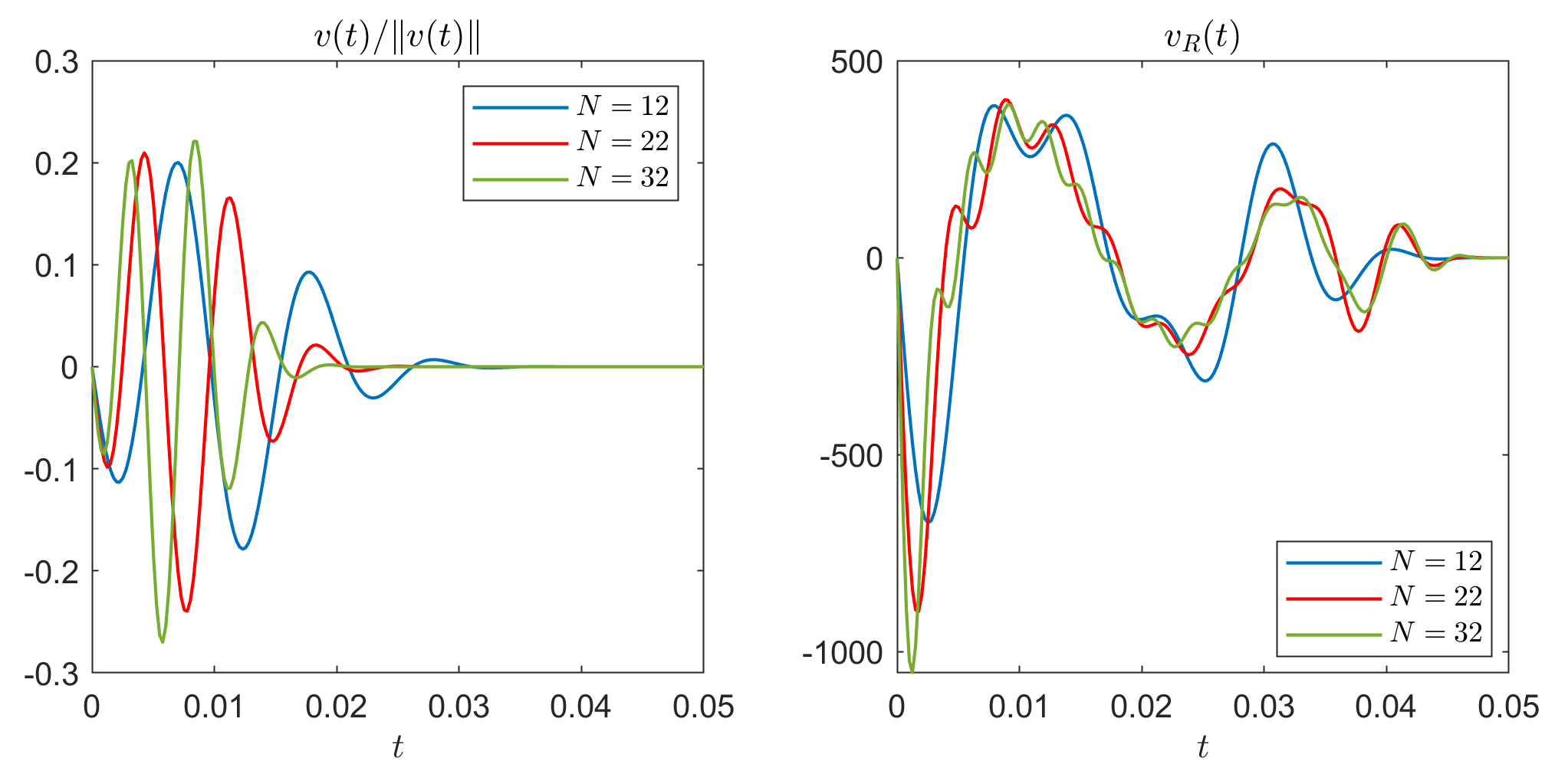}
\end{center}\caption{The exact (left) and the regularized (right) controls of the fourth example for different values of $N$.}\label{fig6}
\end{figure}

\textbf{Example 4:} The Dirichlet-Dirichlet problem ($\alpha=0$) for the heat equation ($K_0=0$) is examined in this example. We consider the sinusoidal initial condition \eqref{initial_sin}, a relatively small final time $T=0.05$ and large $N = 12,\, 22$ and $32,$ a combination of parameters which results in large (in norm) controls. In \autoref{table8} we compare the exact and the regularized $L^2$ norms of the control and the error, namely the solution at $t=T.$ The controls are also presented in \autoref{fig6}. We note that we show only the normalized functions of the exact controls, since their varying order of magnitude does not allow for a combined plot.

\begin{table}[ht]
\begin{center}
 \begin{tabular}{ | c  | c  | c  |c|| c | c|} 
 \hline
 $N$ & $\lVert v(t)\rVert_2$ & $\lVert v_R(t)\rVert_2$ & Rel. change & $\lVert\theta(x,T)\rVert_2$  & $\lVert\theta_R(x,T)\rVert_2$  
\\ \hhline{======}
12 &   714.8 &     53.9  & $-92.5\%$  &  8.324713E-5   & 7.671623E-3 
 \\
   22    &  9436.0 & 52.2  & $-99.5\%$  & 8.912321E-15 & 9.345277E-3 
  \\ 
   32     &  65118.3 &  51.9  &  $-99.9\%$ & 4.993819E-28  &   9.156095E-3 
    \\ \hline
 \end{tabular}
\caption{The fourth example for $K_0 =0$ and $T=0.05.$: The first two columns are the $L^2$ norms of the exact and regularized control. The last two columns are $L^2$ norms of the associated errors, namely the solutions at the final time.}\label{table8}
\end{center}
\end{table}

\section*{Appendix}\label{sec_appendix}
Here we present some qualitative results for the roots of $\Delta(\lambda)$.

First, the definition \eqref{def:Delta} reads
$$\Delta(\lambda)=0 \quad \Rightarrow \quad e^{2i\lambda L} = \dfrac{(1-i \alpha \lambda) \left(1+i \beta \lambda - \beta\frac{K_0}{D_0}\right)}{ \left(1+i \alpha \lambda - \alpha\frac{K_0}{D_0}\right) (1-i \beta \lambda)} e^{\frac{K_0}{D_0}L}.$$
Looking for roots with large modulus we obtain that $ e^{2i\lambda L} \sim e^{\frac{K_0}{D_0}L}, \ \lambda \to \infty,$ resulting in
$$\lambda_n \sim - i\frac{K_0}{2D_0} + \frac{n \pi}{L}, \qquad n\to \infty,$$
which reside outside of the domains $D_\pm$. Hence, possible roots that reside in $D_\pm$ are of finite modulus. Thus, it is possible to deform/choose $\partial D_\pm$ which satisfy the properties of Remark \ref{ref-deform}, and they are the boundaries of domains which contain no roots of $\Delta(\lambda)$, namely no singularities.

In fact, one could prove that there are at most 4 roots of $\Delta(\lambda)$, which do not lie on the line $- i\frac{K_0}{2D_0} + \R$; in particular, they lie on the imaginary axis $\big(\text{Re} \lambda=0\big)$, and they appear in pairs symmetric to the point $- i\frac{K_0}{2D_0} $.

A rigorous proof is out of the scope of this manuscript, but a sketch of a proof looks as follows. Making the change of variables $\lambda=i\left(\frac{y}{L}-\frac{K_0}{2D_0}\right)$, the equation $\Delta(\lambda)=0$ takes the form
$$\frac{\sigma y}{1+\rho y^2}=\tanh y, \qquad \sigma\in \R, \ \ \rho\in\R,$$
which admits solutions only if $y$ is real or imaginary. The latter case, corresponds to the set of roots lying on the line $- i\frac{K_0}{2D_0} + \R$, which is outside of $D_\pm$. Thus, one needs to examine only the case of $y\in\R$.

Without proof we present in \autoref{table9} the number of solutions of $\Delta(\lambda)=0$ with respect to the parameters $\rho$ and $\sigma$. For $\rho>0$ and $0<\sigma<1,$ the number of roots  depends on the sign of the quantity $\sigma-2\sqrt{\rho} \tanh\frac{1}{\sqrt{\rho}}$, approximately.\footnote{The exact quantity is $\sigma-\left(\frac{\rho }{\sqrt{\eta }}+\sqrt{\eta }\right) \tanh \left(\frac{\sqrt{\eta }}{\rho }\right)$, where $\eta=\frac{1}{2} \sigma  (\sigma -\rho )-\rho+\sqrt{\left(\frac{1}{2} \sigma  (\sigma -\rho )-\rho \right)^2-(1-\sigma ) \rho ^2}$.} Using the fact that $|\tanh\frac{1}{\sqrt{\rho}}|<1$, if $\sigma>2\sqrt{\rho} $, then 4 roots occur.

\begin{table}[ht]
\begin{center}
\begin{tabular}{ |c|c||c| } 
\hline
$\rho$ & $\sigma$ & No. of roots  \\ \hhline{===}
\multirow{3}{4em}{$\phantom{=}>0$} & $(-\infty,\, 0] $ & $0$ \\ 
& $(0,\,1)$ & $\{0,\,2,\,4\}$ \\ 
& $[1,\, \infty)$ & $2$ \\ 
\hline
\multirow{3}{4em}{$\phantom{=}=0$} & $(-\infty,\, 0] $ & $0$ \\ 
& $(0,\,1)$ & $2$ \\ 
& $[1,\, \infty)$ & $0$ \\ 
\hline
\multirow{2}{4em}{$\phantom{=}<0$} & $(-\infty,\, 1) $ & $2$ \\ 
& $[1,\,\infty)$ & $0$ \\ 
\hline
\end{tabular}
\caption{The number of roots of $\Delta (\lambda)$ for different values of $\rho$ and $\sigma.$ }\label{table9}
\end{center}
\end{table}

%
For matters of completeness we present the exact values of $\sigma$ and $\rho$ in terms of the parameters of the IBVP:
$$\sigma=\frac{1}{L}\dfrac{\alpha-\beta}{\left(1-\alpha\frac{K_0}{2D_0}\right)\left(\beta\frac{K_0}{2D_0}-1\right)} \qquad \text{and} \qquad \rho =\frac{1}{L^2}\dfrac{\alpha\,\beta}{\left(1-\alpha\frac{K_0}{2D_0}\right)\left(\beta\frac{K_0}{2D_0}-1\right)}.$$

\bibliographystyle{siam}
\bibliography{refs}

\end{document}